\documentclass[12pt,a4paper]{article}
\usepackage[T1]{fontenc}
\usepackage[utf8]{inputenc}
\usepackage{amsmath,amssymb} 
\usepackage{newtxtext, newtxmath}
\usepackage{microtype}
\usepackage[hidelinks]{hyperref}
\usepackage[margin=1.2in]{geometry}
\usepackage{setspace}
\newtheorem{theorem}{Theorem}[section]

\numberwithin{equation}{section}

\begin{document}

\thispagestyle{plain}

\noindent{Revised November 11, 2026,  with an addendum at the end.}\\
\noindent{Journal Information: J. Symbolic Comput. 44  (2009)  960--971.}

\vspace{3mm}

\begin{center}
{\Large\bfseries The $q$-WZ Method for Infinite Series\par}

\vspace{2mm}

William Y.C. Chen\textsuperscript{1}
and Ernest X.W. Xia\textsuperscript{2}\par

Center for Combinatorics, LPMC-TJKLC\\
Nankai University\\
Tianjin 300071, P.R. China\par

\vspace{3mm}

{\small Emails: \textsuperscript{1}chen@nankai.edu.cn,
\textsuperscript{2}xxw@cfc.nankai.edu.cn}
\end{center}

\begin{abstract}
Motivated by the telescoping proofs of two identities of Andrews and
Warnaar, we find that infinite $q$-shifted factorials can be incorporated
into the implementation of the $q$-Zeilberger algorithm in the approach
of Chen, Hou and Mu to prove nonterminating basic hypergeometric series
identities. This observation enables us to extend the $q$-WZ method to
identities on infinite series. We give the $q$-WZ pairs for some
classical identities such as the $q$-Gauss sum, the ${}_6\phi_5$ sum,
the Ramanujan's ${}_1\psi_1$ sum and Bailey's ${}_6\psi_6$ sum.
\end{abstract}


\section{Introduction}
The objective of this paper is to give an extension of the $q$-WZ method (Amdeberhan and Zeilberger, 1998; Mohammed and Zeilberger, 2004; Wilf and Zeilberger, 1990a; Zeilberger, 1990, 1994) to nonterminating basic hypergeometric series identities. We will follow the standard notation on $q$-series (Gasper and Rahman, 2004) and always assume $0<|q|<1$. The $q$-shifted factorials $(a;q)_n$ and $(a;q)_\infty$ are defined by
\begin{eqnarray*}
(a;q)_n & = &\left\{
     \begin{array}{ll}
         1, \qquad & \mbox{if $n=0$}, \\ [2pt]
        (1-a) (1-aq) \cdots (1-aq^{n-1}), \qquad
              & \mbox{if $n\geq 1$},
     \end{array} \right. \\ [3pt]
     (a;q)_{-n}&=&\frac{1}{(aq^{-n};q)_n},\\[2pt]
(a;q)_\infty & = & (1-a) (1-aq)  (1-aq^2) \cdots,\\[2pt]
 (a_1,a_2,\ldots,a_k;q)_n&=&(a_1;q)_n(a_2;q)_n\cdots (a_k;q)_n.
\end{eqnarray*}
An ${}_r\phi_s$ basic hypergeometric series is defined by
\begin{align}\label{phi}
_r\phi_{s}\left[{a_1,a_2,\ldots,a_r \atop
b_1,b_2,\ldots,b_{s}};q,z\right]:=
\sum_{n=0}^\infty\frac{(a_1,a_2,\ldots,a_r;q)_n} {(q,b_1,\ldots,
b_{s};q)_n}\left[(-1)^nq^{n\choose 2}\right]^{1+s-r}z^n,
\end{align}
where $q\neq0$ when $r>s+1$. An ${}_r\psi_s$ bilateral basic hypergeometric series is defined by
\begin{align}\label{psi}
_r\psi_s\left[ {a_1,a_2,\ldots,a_r \atop b_1,b_2,\ldots,b_s};q,z
\right]:=\sum_{n=-\infty}^\infty\frac{(a_1,a_2,\ldots,a_r;q)_n}
{(b_1,b_2,\ldots,b_s;q)_n}\left[(-1)^nq^{n\choose
2}\right]^{s-r}z^n.
\end{align}
It is assumed that $q$, $z$ and the parameters are such that each term of the series is well-defined. We say that an ${}_r\phi_s$ series terminates if only a finite number of terms contribute. Otherwise, we say that the series ${}_r\phi_s$ is nonterminating.

For the ordinary nonterminating hypergeometric identities, Gessel (1995) and Koornwinder (1998) provided computer proofs of Gauss' summation formula and Saalsch\"utz's summation formula by means of a combination of Zeilberger's algorithm and asymptotic estimates. Vidunas (2002) (see also Koepf (2003) and Koornwinder (2005)) presented a method to evaluate ${}_2F_1\bigl({a,b\atop c}\bigm|-1\bigr)$ when $c-a+b$ is an integer. Recently, Chen et al.\ (2008) developed an approach to proving nonterminating basic hypergeometric identities based on the $q$-Zeilberger algorithm (Koornwinder, 1993; Wilf and Zeilberger, 1992). In this paper we will show how to apply the $q$-WZ method to prove nonterminating basic hypergeometric summation formulas by finding the $q$-WZ pairs. We will give some examples including the $q$-Gauss sum, the very-well-poised ${}_6\phi_5$ sum, the Ramanujan's ${}_1\psi_1$ sum and Bailey's very-well-poised series ${}_6\psi_6$ sum (Gasper and Rahman, 2004).

\section{The Andrews--Warnaar identities}
In this paper, we give telescoping proofs of the following two identities on partial theta functions (Andrews, 1981):
\begin{eqnarray}
\bigg(\sum_{n=0}^\infty(-1)^na^nq^{n\choose2
}\bigg)\bigg(\sum_{n=0}^\infty(-1)^nb^nq^{n\choose2
}\bigg)&=&(q,a,b;q)_\infty\sum_{n=0}^\infty\frac{(abq^{n-1};q)_n}
{(q,a,b;q)_n}q^n,
\label{Andrews-Warnaar-1}\\[2pt]
 1+\sum_{n=1}^\infty(-1)^nq^{n\choose
2}(a^n+b^n)&=&(a,b,q;q)_\infty \sum_{n=0}^\infty\frac{(ab/q;q)_{2n}}
{(q,a,b,ab;q)_n}q^n. \label{Andrews-Warnaar-2}
\end{eqnarray}
The identity \eqref{Andrews-Warnaar-2} was first proved by Warnaar (2003). Andrews and Warnaar (2007) derived the identity \eqref{Andrews-Warnaar-1} and used it to prove \eqref{Andrews-Warnaar-2}.

As will be seen, the telescoping proofs suggest that the approach developed by Chen et al.\ (2008) for proving nonterminating basic hypergeometric identities can be extended so that infinite $q$-shifted factorials can be allowed in a $q$-hypergeometric term. This idea immediately leads to an extension of the $q$-WZ method to identities on infinite series.

Note that formula \eqref{Andrews-Warnaar-2} is a generalization of the well-known Jacobi's triple product identity. When $b=q/a$, we get Jacobi's triple product identity
\begin{eqnarray}
\sum_{n=-\infty}^\infty(-1)^na^nq^{n\choose 2}=(a,q/a,q;q)_\infty,
\label{Jacobi-s-triple-product-identity}
\end{eqnarray}
where $|q|<1$ and $a\neq0$.

We now describe how to prove the identities \eqref{Andrews-Warnaar-1} and \eqref{Andrews-Warnaar-2} by the telescoping method. Let us consider \eqref{Andrews-Warnaar-1} first. Put
\[
f(a)=\left(\sum_{n=0}^{\infty}(-1)^nq^{\binom n2}a^n\right)
\left(\sum_{n=0}^{\infty}(-1)^nq^{\binom n2}b^n\right).
\]
Note that the second factor does not contain the parameter $a$. It is easily verified that
\begin{align}\label{R1}
f(a)=(1-a)f(aq)+aqf(aq^2).
\end{align}
We proceed to show that the right hand side of \eqref{Andrews-Warnaar-1} satisfies the same recurrence relation. Of course, we still need to verify the boundary condition. Let
\[
g(a)=\sum_{n=0}^{\infty}D_n(a),\quad\text{where }\; D_n(a)=(q,a,b;q)_\infty
\frac{(abq^{n-1};q)_nq^n}{(q,a,b;q)_n}.
\]
Then it is necessary to show that
\begin{equation} \label{fa} g(a) -(1-a)g(aq) -aqg(aq^2)=0.
\end{equation}
Here comes the key step of finding a telescoping relation for $D_n(a)$. Note that, for any $n\ge0$, we have
\begin{align}\label{g1}
&D_n(a) -(1-a)D_n(aq)
-aqD_n(aq^2)\nonumber\\[2pt]
=&\frac{(abq^n;q)_n(q,a,b;q)_\infty q^n}{(q,a,b;q)_n} \bigg(
\frac{1-abq^{n-1}}{1-abq^{2n-1}}
-\frac{1-a}{1-aq^n}-\frac{aq(1-abq^{2n})}
{(1-aq^{n+1})(1-aq^n)(1-abq^n)} \bigg) \nonumber\\[2pt]
=&\frac{(abq^n;q)_n(q,a,b;q)_\infty q^n}{(q,a,b;q)_n} \bigg(
\frac{a(1-q^n)(1-bq^{n-1})}{(1-aq^n)(1-abq^{2n-1})}-
\frac{aq(1-abq^{2n})} {(1-aq^{n+1})(1-aq^n)(1-abq^n)}
\bigg) \nonumber\\[2pt]
=&z_{n+1}-z_n,
\end{align}
where
\[
z_n=-\frac{(1-q^n)(1-bq^{n-1})(abq^n;q)_n(q,a,b;q)_\infty aq^n}
{(1-aq^n)(1-abq^{2n-1})(q,a,b;q)_n}.
\]
The above relation reveals that the infinite $q$-shifted factorial $(q,a,b;q)_\infty$ can be incorporated into the telescoping relation and this step can be automated by the $q$-Gosper algorithm (Gosper, 1978; Koornwinder, 1993). Moreover, one sees that infinite $q$-shifted factorials can be incorporated into the $q$-Zeilberger algorithm so that the approach of Chen et al.\ (2008) can be extended to terms containing infinite $q$-shifted factorials. In particular, one can make the $q$-WZ method work for nonterminating hypergeometric series.

Now, let us return to the proof of \eqref{Andrews-Warnaar-1}. Clearly, $z_0=0$. It is also easily seen that $\lim_{n\to+\infty}z_n=0$. Summing \eqref{g1} over the nonnegative integers, we obtain the recurrence relation \eqref{fa}. In order to show that $f(a)=g(a)$, we will use the recurrence relation of $f(a)-g(a)$ to reach this goal.

Let $H(a)=f(a)-g(a)$. From the recurrence relations for $f(a)$ and $g(a)$, it follows that $H(a)$ satisfies the recurrence relation
\begin{align}\label{recurrence-relation1}
H(a)=(1-a)H(aq)+aqH(aq^2).
\end{align}
Iterating the above relation yields that
\begin{align}{\label{Recurrence}}
H(a)=A_nH(aq^{n+1})+B_nH(aq^{n+2}),
\end{align}
where $A_n$ and $B_n$ are given by
\[
A_0=(1-a),\quad B_0=aq,\quad A_1=(1-a)(1-aq)+aq,\quad B_1=(1-a)aq^2,
\]
and
\[
A_{n+1}=(1-aq^{n+1})A_n+aq^{n+1}A_{n-1},\quad
B_{n+1}=aq^{n+2}A_n,\quad n\ge1.
\]
Hence we have
\[A_{n+1}-A_n=-aq^{n+1}(A_n-A_{n-1}),\]
which implies that
\begin{align*}
|A_{n+1}-A_n|&=|(-1)^na^nq^{\binom{n+2}{2}-1}|\,|(A_1-A_0)|\\
&\le |a^nq^{\binom{n+2}{2}-1}|\,(|A_1|+|A_0|).
\end{align*}
So, for fixed $a$ and $|q|<1$, $\lim_{n\to+\infty}A_n$ exists. Since $B_{n+1}=aq^{n+2}A_n$, $\lim_{n\to+\infty}B_n$ also exists. Again, by relation \eqref{Recurrence}, we find
\[
H(a)=H(0)\left(\lim_{n\to+\infty}A_n+\lim_{n\to+\infty}B_n\right).
\]
It remains to show that $H(0)=0$, that is,
\begin{equation}\label{h0}
 \sum_{n=0}^\infty(-1)^nb^nq^{n\choose 2}
=(q,b;q)_\infty\sum_{n=0}^\infty\frac{q^n}{(q,b;q)_n}.
 \end{equation}
We can use the telescoping method to prove \eqref{h0}. Let
\[
G(b)=\sum_{n=0}^{\infty}(-1)^nb^nq^{\binom n2}
-(q,b;q)_\infty\sum_{n=0}^{\infty}\frac{q^n}{(q,b;q)_n}.
\]
It is easy to check that
\[G(b)=(1-b)G(bq)+bqG(bq^2).\]
We aim to show that $G(b)=0$. Since $G(b)$ satisfies the same recurrence relation as $H(a)$, it suffices to confirm $G(0)=0$, that is,
\[(q;q)_\infty\sum_{n=0}^{\infty}\frac{q^n}{(q;q)_n}=1,\]
which is a special case of Euler's identity (Gasper and Rahman, 2004, P.~354)
\[
\sum_{n=0}^{\infty}\frac{z^n}{(q;q)_n}=\frac1{(z;q)_\infty},\quad |z|<1.
\]
Indeed, relation \eqref{h0} is a limiting case of Heine's transformation of ${}_2\phi_1$. For completeness, we give a proof based on Euler's identities:
\begin{eqnarray}\label{identity-1}
(q,b;q)_\infty\sum_{n=0}^\infty\frac{q^n}{(q,b;q)_n} &=&
(q;q)_\infty\sum_{m=0}^\infty
\frac{q^m}{(q;q)_m}\sum_{n=0}^\infty\frac{q^{n\choose
2}(-bq^m)^n}{(q;q)_n}\nonumber\\[2pt]
&=&(q;q)_\infty\sum_{n=0}^\infty\frac{(-1)^nb^nq^{n\choose
2}}{(q;q)_n}\sum_{m=0}^\infty\frac{(q^{n+1})^m}{(q;q)_m}\nonumber
\\[2pt]
&=& \sum_{n=0}^\infty(-1)^nb^nq^{n\choose 2}.
\end{eqnarray}
Thus, we have verified that $H(a)=0$. This completes the proof.

We remark that once the recurrence relation \eqref{recurrence-relation1} is derived, one can also use the theorem of Chen et al.\ (2008, Theorem~3.1) to prove the existence of the limits of $A_n$ and $B_n$.

We next present a telescoping proof of \eqref{Andrews-Warnaar-2}. Let
\[f(a)=1+\sum_{n=1}^{\infty}(-1)^nq^{\binom n2}(a^n+b^n).\]
It is easily seen that
\begin{align}\label{R2}
(1+aq)f(a)-(1-a^2q)f(aq)-(aq+a^2q)f(aq^2)=(q-1)a.
\end{align}
Let
\[
g(a)=\sum_{n=0}^{\infty}D_n(a),\quad\text{where }D_n(a)=(q,a,b;q)_\infty
\frac{(ab/q;q)_{2n}q^n}{(q,a,b,ab;q)_n}.
\]
It will be shown that
\begin{align}\label{ga}
(1+aq)g(a)-(1-a^2q)g(aq)-(aq+a^2q)g(aq^2)=(q-1)a.
\end{align}
Since
\begin{align*}
&\frac{q^n-abq^{n-1}}{1-abq^{2n-1}}
-\frac{(1-a^2q)(1-ab)q^n}{(1+aq)(1-aq^n)
(1-abq^n)}\\[2pt]
&\qquad-\frac{(a^2q+aq)(1-abq^{2n})(1-abq)q^n}
{(1+aq)(1-aq^n)(1-aq^{n+1})
(1-abq^n)(1-abq^{n+1})}\\[2pt]
&\ =\frac{(1-abq^{2n})(-1+q+abq^{n+1}+a^2bq^{n+2}
-aq^{n+2}-q^{n+1}-a^2bq^{2n+2} +a^2bq^{2n+3})a}{
(1-aq^{n+1})(1-aq^n)(1-abq^{n})(1+aq)(1-abq^{n+1})}\\[2pt]
&\ -\frac{
(-1+q+abq^n+a^2bq^{n+1}-aq^{n+1}-q^n-a^2bq^{2n}+a^2bq^{2n+1})
a(1-bq^{n-1})(1-q^n)} {(1-aq^n)(1-abq^{2n-1})(1+aq)(1-abq^n)},
\end{align*}
multiplying both sides by
\[\frac{(ab;q)_{2n}(a,b,q;q)_\infty}{(q,a,b,ab;q)_n},\]
we deduce that
\begin{align}\label{g2}
&D_n(a) -\frac{(1-a^2q)}{1+aq} D_n(aq)
-\frac{(a^2q+aq)}{1+aq}D_n(aq^2) =z_{n+1}-z_{n},
\end{align}
where
\begin{align*}
z_n=&\frac{(-1+q-aq^{n+1}-q^n+a^2bq^{n+1}+abq^n-a^2bq^{2n}
+a^2bq^{2n+1})a}{(1-aq^n)(1-abq^{2n-1})}\\[2pt]
 &\qquad\qquad\times\frac{(1-bq^{n-1})(1-q^n)(ab;q)_{2n}(q,a,b;q)_\infty}
{(1+aq)(1-abq^n)(ab;q)_n(q,a,b;q)_n}.
\end{align*}
Clearly, $z_0=0$ and $\lim_{n\to+\infty}z_n=\frac{(q-1)a}{1+aq}$. Summing \eqref{g2} over nonnegative integers, we obtain the recurrence relation \eqref{ga}.

Let $H(a)=f(a)-g(a)$. Then $H(a)$ satisfies the following recurrence relation
\begin{align}\label{recurrence-relation2}
H(a)=\frac{1-a^2q}{1+aq}H(aq)+ \frac{aq+a^2q}{1+aq}H(aq^2).
\end{align}
By iteration, we obtain
\begin{align}{\label{Recurrence2}}
H(a)=A_nH(aq^{n+1})+B_nH(aq^{n+2}),
\end{align}
where $A_n$ and $B_n$ are given by
\begin{align*}
A_0&=\frac{1-a^2q}{1+aq},& A_1&=\frac{1+a^3q^3}{1+aq^2},\\
B_0&=\frac{aq+a^2q}{1+aq},& B_1&=\frac{aq^2(1-a^2q)}{1+aq^2},
\end{align*}
and for $n\ge1$,
\begin{eqnarray}
 A_{n+1}& = & \frac{1-a^2q^{2n+3}}{1+aq^{n+2}}A_{n}
+\frac{aq^{n+1}+a^2q^{2n+1}}{1+aq^{n+1}}A_{n-1},\label{Formula-1} \\[2pt]
B_{n+1}& =& \frac{aq^{n+2}+a^2q^{2n+3}}{1+aq^{n+2}}A_n.
\label{Formula-2}
\end{eqnarray}
Based on the above recurrence relations, one can deduce that both $\lim_{n\to+\infty}A_n$ and $\lim_{n\to+\infty}B_n$ exist. We note that Zeilberger (2008) has shown that
\[A_n=\frac{1+(-1)^{n+1}a^{n+2}q^{\binom{n+2}{2}}}{1+aq^{n+1}}\]
and
\[B_n=\frac{aq^{n+1}\left(1+(-1)^na^{n+1}q^{\binom{n+1}{2}}\right)}{1+aq^{n+1}}.\]
Now we see that the limits $\lim_{n\to+\infty}A_n$ and $\lim_{n\to+\infty}B_n$ exist. By relation \eqref{Recurrence2}, we deduce that
\[H(a)=H(0)\left(\lim_{n\to+\infty}A_n+\lim_{n\to+\infty}B_n\right).\]
The identity \eqref{identity-1} implies that $f(0)=g(0)$. So we have $H(a)=0$. This completes the proof.

We also note that once the recurrence relation \eqref{recurrence-relation2} is established, one may assume that $|a|<1$ and may use the theorem in Chen et al.\ (2008, Theorem~3.1) to prove the existence of the limits of $A_n$ and $B_n$. Moreover, we may drop the assumption $|a|<1$ by analytic continuation.

\section{The $q$-WZ pairs for infinite series}
Our approach to the $q$-WZ method for infinite series can be described as follows. The key step is to construct $q$-WZ pairs for infinite sums. Suppose that we aim to prove an identity of the form:
\begin{align}\label{equation1}
\sum_{k=N_0}^\infty F_k(a_1,a_2,\ldots, a_t)=R(a_1,a_2,\ldots,a_t),
\end{align}
where $t$ is a positive integer, and the sum is either a unilateral or bilateral basic hypergeometric series, namely, $N_0=0$ or $N_0=-\infty$, $R(a_1,a_2,\ldots,a_t)$ is either zero or a quotient of two products of infinite $q$-shifted factorials.

First, we set some parameters, say, $a_1,\ldots,a_p$, $(1\le p\le t)$ to $a_1q^n,\ldots,a_pq^n$, so that we get
\begin{align}\label{-Set-parameters-}
\sum_{k=N_0}^\infty F_k(a_1q^n,\ldots,
a_pq^n,a_{p+1},\ldots,a_t)=R(a_1q^n,\ldots,
a_pq^n,a_{p+1},\ldots,a_t).
\end{align}
If $R(a_1q^n,\ldots,a_pq^n,a_{p+1},\ldots,a_t)\neq0$, set
\[
F(n,k)=\frac{F_k(a_1q^n,\ldots,a_pq^n,a_{p+1},\ldots,a_t)}
{R(a_1q^n,\ldots,a_pq^n,a_{p+1},\ldots,a_t)}.
\]
If \(R\) is identically zero, set  
\[F(n,k)=F_k(a_1q^n,\ldots,a_pq^n,a_{p+1},\ldots,a_t).\]
Our goal is to show that
\begin{align}\label{Constant}
\sum_{k=N_0}^\infty F(n,k)= {\rm constant}, \ \ \ \ \ \
n=0,1,2,\ldots.
\end{align}
It should be noted that the above constant is meant to be 
independent of $n$ but may depend on the parameters,
while in all the examples the constants turn out to 
be independent of the parameters. 

Let us recall the boundary and limit conditions for the $q$-WZ-method. 
Let $f(n)$ denote the left hand side of \eqref{Constant}, i.e.,
\[f(n)=\sum_{k=N_0}^{\infty}F(n,k)\]
and we aim to show that
\[f(n)=\mathrm{constant}\]
for every nonnegative integer $n$. To this end, it suffices to show that $f(n+1)-f(n)=0$ for every nonnegative integer $n$. This can be done by finding $G(n,k)$ such that
\begin{align}\label{WZPair}
F(n+1,k)-F(n,k)=G(n,k+1)-G(n,k).
\end{align}
A pair of functions $(F(n,k),G(n,k))$ satisfying \eqref{WZPair} is called a $q$-WZ pair. Once a $q$-WZ pair is found, one can check the boundary and limit conditions to ensure that $f(n)$ equals the claimed constant. Here are the conditions:
\begin{list}{}{\setlength{\leftmargin}{2.8em}\setlength{\labelwidth}{2.5em}\setlength{\itemsep}{2pt}\setlength{\parsep}{0pt}}
\item[(C1)] For each integer $n\ge0$, $\lim_{k\to\pm\infty}G(n,k)=0$.
\item[(C2)] For each integer $k$, the limit
\begin{align}\label{limit}
f_k=\lim_{n\rightarrow \infty}F(n,k)
\end{align}
exists and is finite.
\item[(C3)] We have $\lim_{L\to\infty}\sum_{n\ge0}G(n,-L)=0$.
\end{list}
The WZ method can be formally stated as follows.
\begin{theorem}[Wilf and Zeilberger (1990b)]
Assume that $(F(n,k),G(n,k))$ is a WZ pair. If {\rm(C1)} holds  and $\sum_k|F(n,k)|<\infty$ for every $n\ge0$, then we have
\begin{align}\label{Con1}
\sum_{k}F(n,k)={\rm constant},\ \ \ \ \ \ \ n=0,1,2,\ldots.
\end{align}
If {\rm(C2)} and {\rm(C3)} hold, then we have the companion identity
\begin{align}\label{Identity}
\sum_{n=0}^\infty G(n,k)=\sum_{j\leq k-1}(f_j-F(0,j)),
\end{align}
where $f_j$ is defined by \eqref{limit}.
\end{theorem}
We now explain how to compute the desired $q$-WZ pair for the identity \eqref{equation1}. 
In fact, it can be obtained by applying the $q$-Gosper algorithm to $F(n+1,k)-F(n,k)$. 
It should be noted that $F(n+1,k)-F(n,k)$ is a $q$-hypergeometric term with respect to
$q^k$, even if $F(n,k)$ contains infinite $q$-shifted factorials such as $(aq^n;q)_\infty$. 
Obviously, $F(n+1,k)-F(n,k)$ is a $q$-hypergeometric term when $R(a_1,\ldots,a_t)=0$. Assume that $R(a_1,\ldots,a_t)\neq0$. Let
\begin{align*}
 M_1&=\frac{R(a_1q^{n+1},\ldots,
a_pq^{n+1},a_{p+1},\ldots,a_t)}{R(a_1q^n,\ldots,
a_pq^n,a_{p+1},\ldots,a_t)},\\[2pt]
M_2 &=\frac{F_{k+1}(a_1q^{n+1},\ldots,
a_pq^{n+1},a_{p+1},\ldots,a_t)}{F_k(a_1q^{n+1},\ldots,
a_pq^{n+1},a_{p+1},\ldots,a_t)},\\[2pt]
M_3&=\frac{F_{k+1}(a_1q^n,\ldots,
a_pq^n,a_{p+1},\ldots,a_t)}{F_k(a_1q^{n+1},\ldots,
a_pq^{n+1},a_{p+1},\ldots,a_t)},\\[2pt]
M_4&=\frac{F_k(a_1q^n,\ldots,
a_pq^n,a_{p+1},\ldots,a_t)}{F_k(a_1q^{n+1},\ldots,
a_pq^{n+1},a_{p+1},\ldots,a_t)}.
 \end{align*}
Since $M_1$ is a rational function in $q^n$ and is independent of $k$, $M_2,M_3,M_4$ are all rational functions in $q^k$. Observe that
\begin{equation}\label{fm}
 \frac{F(n+1,k+1)-F(n,k+1)}{F(n+1,k)-F(n,k)}=
\frac{M_2-M_1M_3}{1-M_1M_4}
\end{equation}
is a rational function in $q^k$, i.e., $F(n+1,k)-F(n,k)$ is a $q$-hypergeometric term with respect to $q^k$. It is necessary to mention that even if $F(n,k)$ contains infinite $q$-shifted factorials of the form $(aq^n;q)_\infty$, the quotient \eqref{fm} no longer contains the $q$-shifted factorial $(aq^n;q)_\infty$ and it is still a rational function in $q^k$. Consequently, we can employ the $q$-Gosper algorithm to determine whether $G(n,k)$ exists. Nevertheless, it is also necessary to note that $G(n,k)$ contains infinite $q$-shifted factorials if $F(n,k)$ does.

There is another way to look at the above procedure. Suppose that $F(n,k)$ contains an infinite $q$-shifted factorial $(a;q)_\infty$, where $a$ is a chosen parameter for the substitution $a\to aq^n$. If we set $G'(n,k)=R(aq^n)G(n,k)$, then Eq.~\eqref{WZPair} becomes
\[F(n+1,k)R(aq^n)-F(n,k)R(aq^n)=G'(n,k+1)-G'(n,k).\]
It is evident that the infinite $q$-shifted factorial $(aq^n;q)_\infty$ will disappear in the above equation, and one can use the $q$-Gosper algorithm to find a $q$-WZ pair if it exists.

In the unilateral examples below, $F(n,k)$ and $G(n,k)$ are extended by zero for $k<0$. For fixed generic parameters in the stated convergence domains, $|F(n,k)|\le C\rho^k$ for $k\ge0$ and all sufficiently large $n$, with $C>0$ and $0<\rho<1$; thus dominated convergence applies.

We now take the $q$-binomial theorem (Gasper and Rahman, 2004, P.~354) as an example to explain the above steps:
\begin{align}\label{q-binomial-theorem}
\sum_{k=0}^\infty\frac{(a;q)_k}{(q;q)_k}z^k=
\frac{(az;q)_\infty}{(z;q)_\infty},\ \ \ \ \ |z|<1.
\end{align}
In this case, we have
\[F_k(a)=\frac{(a;q)_k}{(q;q)_k}z^k,\qquad R(a)=\frac{(az;q)_\infty}{(z;q)_\infty}.\]
We choose the parameter $a$, and substitute $a$ with $aq^n$. Then we set
\[F(n,k)=\frac{F_k(aq^n)}{R(aq^n)}=\frac{(aq^n;q)_k(z;q)_\infty}{(q;q)_k(azq^n;q)_\infty}z^k.\]
In order to find $G(n,k)$ such that \eqref{WZPair} holds, it is easily checked that $F(n+1,k)-F(n,k)$ is a $q$-hypergeometric term. By examining the $q$-Gosper algorithm, one sees that it is capable of dealing with the input $F(n+1,k)-F(n,k)$, or we can set
\[G'(n,k)=R(aq^n)G(n,k)\]
and find a solution of the equation
\begin{equation}\label{New-WZ-Equation}
(1-azq^n)\frac{(aq^{n+1};q)_k}{(q;q)_k}z^k-\frac{(aq^n;q)_k}{(q;q)_k}z^k
= G'(n,k+1)-G'(n,k).
\end{equation}
Finally, we obtain the $q$-WZ pair
\begin{align*}
F(n,k)&=\frac{(aq^n;q)_k(z;q)_\infty}
{(q;q)_k(azq^n;q)_\infty}z^k,\\[2pt]
G(n,k)&=-\frac{(aq^n;q)_k(z;q)_\infty(a-aq^k)}
{(q;q)_k(azq^n;q)_\infty(1-aq^n)}q^nz^k.
\end{align*}
If $|z|<1$, it is easy to see that $F(n,k)$ and $G(n,k)$ satisfy conditions (C1), (C2) and (C3). By \eqref{Con1},
\[\sum_{k=-\infty}^{\infty}F(n,k)=\sum_{k=0}^{\infty}F(n,k)=\mathrm{constant},\quad n=0,1,2,\ldots.\]

By dominated convergence and Euler's identity, letting $n\to\infty$ shows that the constant equals $1$. Setting $n=0$, we have
\[\sum_{k=0}^{\infty}F(0,k)=\mathrm{constant}=1.\]
By \eqref{Identity}, we get the companion identity of \eqref{q-binomial-theorem}
\begin{align*}
\sum_{j=0}^k\frac{(a;q)_j}{(q;q)_j}z^j=
(az;q)_\infty\sum_{j=0}^k\frac{z^j}{(q;q)_j}+
\frac{az^{k+1}(a;q)_{k+1}}{(q;q)_{k}}\sum_{n=0}^\infty
\frac{(az;q)_n(aq^{k+1};q)_n}{(a;q)_{n+1}}q^n.
\end{align*}
We remark that our algorithm depends on the choice of parameters. For a given choice of parameters, it is not guaranteed that one can find a $q$-WZ pair. Nevertheless, this approach applies to many classical identities.

We now give a few more examples.

\medskip\noindent\textbf{Example 3.1.} The $q$-Gauss sum (Gasper and Rahman, 2004, P.~354):
\begin{align}\label{Gauss}
\sum_{k=0}^\infty
\frac{(a,b;q)_k}{(q,c;q)_k}\left(\frac{c}{ab}\right)^k
=\frac{(c/a,c/b;q)_\infty} {(c,c/ab;q)_\infty},\ \ \ \ |c/ab|<1.
\end{align}
By computation we get the $q$-WZ pair
\begin{align*}
F(n,k)&=\frac{(b,aq^n;q)_k(c/ab,cq^n;q)_\infty}
{(q,cq^n;q)_k(c/a,cq^n/b;q)_\infty}
\left(\frac{c}{ab}\right)^k,\\[2pt]
G(n,k)&=-\frac{(a-aq^k)(b,aq^n;q)_k(c/ab,cq^n;q)_\infty}
{(1-aq^n)(q,cq^n;q)_k(c/a, cq^n/b;q)_\infty}
\left(\frac{c}{ab}\right)^kq^n.
\end{align*}
If $|c/ab|<1$, it is easy to verify that the two functions $(F(n,k),G(n,k))$ satisfy relation \eqref{WZPair} and conditions (C1), (C2) and (C3). By \eqref{Con1}, we have
\[\sum_{k=-\infty}^{\infty}F(n,k)=\sum_{k=0}^{\infty}F(n,k)=\mathrm{constant},\quad n=0,1,2,\ldots.\]
By dominated convergence and \eqref{q-binomial-theorem}, the constant equals $\sum_{j=0}^{\infty}f_j=1$, where $f_j$ is displayed below. Thus we have
\[\sum_{k=0}^{\infty}F(0,k)=\mathrm{constant}=1.\]
After simplification, we obtain the identity \eqref{Gauss}.

By \eqref{Identity}, we obtain the companion identity of \eqref{Gauss}
\begin{align}\label{Gauss-c}
-\sum_{n=0}^\infty \frac{(a-aq^k)(b,aq^n;q)_k(c/ab,cq^n;q)_\infty}
{(1-aq^n)(q,cq^n;q)_k(c/a, cq^n/b;q)_\infty}
\left(\frac{c}{ab}\right)^kq^n =\sum_{j\leq k-1}(f_j-F(0,j)),
\end{align}
where
\[f_j=\lim_{n\to\infty}F(n,j)=\frac{(b;q)_j(c/ab;q)_\infty}{(q;q)_j(c/a;q)_\infty}\left(\frac{c}{ab}\right)^j,\]
which can be restated as
\begin{align*}
\sum_{j=0}^k
\frac{(a,b;q)_j}{(q,c;q)_j}\left(\frac{c}{ab}\right)^j=&
\frac{(c/b;q)_\infty}{(c;q)_\infty}
\sum_{j=0}^k\frac{(b;q)_j}{(q;q)_j}
\left(\frac{c}{ab}\right)^j\\[2pt]
& \qquad +\frac{(a,b;q)_{k+1}c^{k+1}}{(q;q)_k(c;q)_{k+1}a^kb^{k+1}}
\sum_{n=0}^\infty\frac{(aq^{k+1},c/b;q)_n}
{(a;q)_{n+1}(cq^{k+1};q)_n}q^n.
\end{align*}
\medskip\noindent\textbf{Example 3.2.} The sum of a very-well-poised ${}_6\phi_5$ series (Gasper and Rahman, 2004, P.~356):
\begin{align}\label{65}
\sum_{k=0}^\infty &
\frac{(1-aq^{2k})(a,b,c,d;q)_k}{(1-a)(q,aq/b,aq/c,aq/d;q)_k}
\left(\frac{aq}{bcd}\right)^k \nonumber\\[2pt]
&\qquad\qquad =\frac{(aq,aq/bc,aq/bd,aq/cd;q)_\infty}
{(aq/b,aq/c,aq/d,aq/bcd;q)_\infty},\ \ \ \ \ \ \ |aq/bcd|<1.
\end{align}
By computation we get the following $q$-WZ pair:
\begin{align*}
F(n,k)=&\frac{(1-aq^{n+2k})(c,d, aq^n,bq^n;q)_k}
{(1-aq^n)(q,aq/b,aq^{n+1}/c,aq^{n+1}/d;q)_k}\\[2pt]
&\qquad \times\frac{(aq/b,aq/bcd,aq^{n+1}/c,aq^{n+1}/d;q)_\infty}
{(aq/bc,aq/bd,aq^{n+1},aq^{n+1}/cd;q)_\infty}
\left(\frac{aq}{bcd}\right)^k,\\[2pt]
G(n,k)=&\frac{(c,d;q)_k(a/b,a/bcd;q)_\infty}
{(q,a/b;q)_k(aq^n,aq^n/cd;q)_\infty}\\[2pt]
&\qquad\times \frac{(aq^n,bq^n;q)_k(aq^n/c,aq^n/d;q)_\infty}
{(aq^n/c,aq^n/d;q)_k(a/bd,a/bc;q)_\infty
}\\[2pt]
&\qquad\times \frac{(a-bc)(a-bd)(aq^n-cd)(1-q^k)}{(a-bcd)(bq^n-1)
(aq^{n+k}-c)(aq^{n+k}-d)} \left(\frac{aq}{bcd}\right)^kq^n.
\end{align*}
It is easily seen that $F(n,k)$ and $G(n,k)$ satisfy conditions (C1), (C2) and (C3). Therefore, by \eqref{Con1}, we have the result that $\sum_{k=0}^{\infty}F(n,k)$ is a constant. By dominated convergence and \eqref{Gauss}, the constant equals $\sum_{k=0}^{\infty}f_k=1$, where $f_k$ is displayed below. Thus we have
\[\sum_{k=0}^{\infty}F(0,k)=\mathrm{constant}=1,\]
which is nothing but \eqref{65}. Since
\[f_k=\frac{(c,d;q)_k(aq/b,aq/bcd;q)_\infty}{(q,aq/b;q)_k(aq/bc,aq/bd;q)_\infty}\left(\frac{aq}{bcd}\right)^k\]
and
\[
F(0,j)=\frac{(1-aq^{2j})(a,b,c,d;q)_j(aq/b,aq/c,aq/d,aq/bcd;q)_\infty}
{(1-a)(q,aq/b,aq/c,aq/d;q)_j(aq,aq/bc,aq/bd,aq/cd;q)_\infty}\left(\frac{aq}{bcd}\right)^j,
\]
by \eqref{Identity}, we obtain the companion identity
\begin{align*}
\sum_{j=0}^k&\frac{(1-aq^{2j})(a,b,c,d;q)_j }
{(1-a)(q,aq/b,aq/c,aq/d;q)_j}
\left(\frac{aq}{bcd}\right)^{j}\\[2pt]
&\qquad\qquad=\frac{(aq,aq/cd;q)_\infty}{(aq/c,aq/d;q)_\infty}
\sum_{j=0}^k\frac{(c,d;q)_j}
{(q,aq/b;q)_j}\left(\frac{aq}{bcd}\right)^j\\[2pt]
&\qquad\qquad+\frac{b(aq;q)_k(b,c,d;q)_{k+1}}
{(q,aq/b;q)_k(aq/c,aq/d;q)_{k+1}}
\left(\frac{aq}{bcd}\right)^{k+1}\\[2pt]
&\qquad\qquad\times\sum_{n=0}^\infty
\frac{(aq/cd;q)_n(aq^{k+1},bq^{k+1};q)_n}{(b;q)_{n+1}
(aq^{k+2}/c,aq^{k+2}/d;q)_{n}}q^n.
\end{align*}
\medskip\noindent\textbf{Example 3.3.} The Ramanujan's ${}_1\psi_1$ sum (Gasper and Rahman, 2004, P.~357)
\begin{align}\label{11}
_1\psi_1(a;b;q,z)=
\frac{(q,b/a,az,q/az;q)_\infty}{(b,q/a,z,b/az;q)_\infty}, \ \ \
|b/a|<|z|<1.
\end{align}
In this case, we find that
\begin{align*}
F(n,k)&=\frac{(aq^n;q)_k(z,b/az,bq^n,q^{1-n}/a;q)_\infty}
{(bq^n;q)_k(q,b/a,azq^n,q^{1-n}/az;q)_\infty}z^k,\\[2pt]
G(n,k)&=\frac{(z,b/az,bq^n,q^{-n}/a;q)_\infty(aq^n;q)_k(1-azq^n)}
{(q,b/a,azq^n,q^{-n}/az;q)_\infty(bq^n;q)_k(z-azq^n)}z^k.
\end{align*}
If $|b/a|<|z|<1$, utilizing the following relation
\begin{align}{\label{Trans}}
(a;q)_{-n}=\frac{(-q/a)^nq^{n\choose 2}}{(q/a;q)_n}, \ \ \ \
n=0,1,2,\ldots,
\end{align}
we can verify that $G(n,k)$ satisfies condition (C1). It follows that
\begin{align} \label{Ra1}
\sum_{k=-\infty}^\infty F(n,k)= {\rm constant}, \ \ \ \ \
n=0,1,2,\ldots.
\end{align}
When \(b\) tends to \(q\), for $k < 0$, we see that the terms vanish, so the \(q\)-binomial theorem  \eqref{q-binomial-theorem} gives \(\lim_{b \to q}\sum_{k=-\infty}^{\infty} F(0,k)=1\).
 Applying \eqref{Ra1} with parameters $(aq^{-m},q,z)$ and $n=m$ proves \eqref{11} for $b=q^{m+1}$, for all sufficiently large integers $m$. For fixed generic $a,z$, both sides of \eqref{11} are holomorphic in $b$ near $0$, with locally uniform convergence of the series. The identity theorem, followed by meromorphic continuation, establishes \eqref{11} throughout its convergence domain.

\medskip\noindent\textbf{Example 3.4.} The sum of a very-well-poised ${}_6\psi_6$ series (Gasper and Rahman, 2004, P.~357):
\begin{align}\label{66}
\sum_{k=-\infty}^\infty&\frac{(1-aq^{2k})(b,c,d,e;q)_k}
{(1-a)(aq/b,aq/c,aq/d,aq/e;q)_k}
\left(\frac{a^2q}{bcde}\right)^k \nonumber\\[2pt]
&\qquad\qquad
=\frac{(aq,aq/bc,aq/bd,aq/be,aq/cd,aq/ce,aq/de,q,q/a;q)_\infty}
{(aq/b,aq/c,aq/d,aq/e,q/b,q/c,q/d,q/e,a^2q/bcde;q)_\infty},
\end{align}
where $\left|a^2q/(bcde)\right|<1$.

We obtain the following $q$-WZ pair:
\begin{align*}
F(n,k)=&\frac{(1-aq^{n+2k})(d,e,bq^n,cq^n;q)_k (aq/b,aq/c;q)_\infty}
{(1-aq^n)(aq/b,aq/c,aq^{n+1}/d,aq^{n+1}/e;q)_k
(aq/bd,aq/be;q)_\infty}\\[2pt]
& \qquad\times \frac{(q/d,q/e,a^2q/bcde,aq^{n+1}/d,
aq^{n+1}/e,q^{1-n}/b,q^{1-n}/c;q )_\infty}
{(q,aq/cd,aq/ce,aq^{n+1},aq^{n+1}/de,
q^{1-n}/a,aq^{1-n}/bc;q)_\infty} \left(\frac{a^2q}
{bcde}\right)^k,\\[2pt]
G(n,k)=&\frac{(d,e,bq^n,cq^n;q)_k(a/b,a/c,1/e,a^2/bcde,
1/d;q)_\infty}{(a/b,a/c,aq^n/d,aq^n/e;q)_k
(q,a/bd,a/be,a/cd,a/ce;q)_\infty }\\[2pt]
&\qquad\times \frac{(aq^n/d,aq^n/e,q^{-n}/b ,q^{-n}/c;q)_\infty
(-1+aq^n)} {
(aq^n,aq^n/de,aq^{-n}/bc,q^{-n}/a;q)_\infty(1-bq^n) (1-cq^n)}\\[2pt]
&\qquad\times\frac{(a-bd)(a-be)(a-cd)(a-ce)(aq^n-de)q^n}
{(aq^{n+k}-d)(aq^{n+k}-e)(a-ad)(1-e)(a^2-bcde)} \left(\frac{a^2q}
{bcde}\right)^k.
\end{align*}
Since $|a^2q/bcde|<1$, from the identity \eqref{Trans} it follows that $G(n,k)$ satisfies condition (C1). By \eqref{Con1}, we find
\begin{align} \label{Ra2}
\sum_{k=-\infty}^\infty F(n,k)= {\rm constant}, \ \ \ \ \
n=0,1,2,\ldots.
\end{align}

When $b$ tends to $a$, we see that for $k<0$,
\(F(0,k)=0\). For \(k\ge0\), \(F(0,k)\) becomes
\[\frac{(1-aq^{2k})(a,c,d,e;q)_k}
{(1-a)( q,aq/c,aq/d,aq/e;q)_k
}\frac{(aq/c,aq/d,aq/e,aq/cde;q)_\infty}
{(aq,aq/cd,aq/ce,aq/de;q)_\infty}\left(\frac{aq}{cde}\right)^k.\] 
The \({}_6\phi_5\) summation formula   \eqref{65} gives
\begin{align}
\begin{aligned}
\lim_{b \to a}\sum_{k=-\infty}^\infty F(0,k) =
&\sum_{k=0}^\infty\frac{(1-aq^{2k})(a,c,d,e;q)_k}
{(1-a)( q,aq/c,aq/d,aq/e;q)_k
}\\[2pt]
& \quad \times \frac{(aq/c,aq/d,aq/e,aq/cde;q)_\infty}
{(aq,aq/cd,aq/ce,aq/de;q)_\infty}\left(\frac{aq}{cde}\right)^k=1.
\label{eq:first}
\end{aligned}
\end{align}

By the symmetry of $b,c,d,e$, \eqref{Ra2} also applies to simultaneous shifts of $a,c,d$. Let $n=m$. Starting with $(aq^{-m},aq^{-m},cq^{-m},dq^{-m},e)$ and using~\eqref{eq:first} proves \eqref{66} for $b=aq^{-m}$ for all sufficiently large integers $m$.  For fixed generic $a,c,d,e$, both sides of \eqref{66} are holomorphic in $t=1/b$ near $0$. 
Indeed, after the cancellation of powers of $t$, choose $r>0$
sufficiently small that all denominators are nonzero on $|t|\le r$
and $\left|a^2q/(cde)\right|r<1$. The terms with indices
$\pm k$, $k\ge1$, are bounded on this disk by
$
C_r\left(\left|\frac{a^2q}{cde}\right|r\right)^k,
$
where $C_r>0$ may depend on $r$ and the fixed parameters,
but is independent of $k$ and $t$. Hence both tails converge
uniformly on $|t|\le r$.
 Since $t=q^m/a\to0$, the identity theorem and meromorphic continuation establish \eqref{66} throughout its convergence domain. 

\section*{Acknowledgments}
The authors would like to thank George E. Andrews, Q.H. Hou and Doron Zeilberger for valuable comments. This work was supported by the 973 Project, the PCSIRT Project of the Ministry of Education, the Ministry of Science and Technology, and the National Science Foundation of China.

\section*{References}
\begingroup\fontsize{11pt}{9.5}\selectfont
\begin{list}{}{\setlength{\leftmargin}{1em}\setlength{\itemindent}{-1em}\setlength{\itemsep}{1pt}\setlength{\parsep}{2pt}\setstretch{1.2}}

\item Amdeberhan, T., Zeilberger, D., 1998. $q$-Apery irrationality proofs by $q$-WZ pairs. Adv. Appl. Math. 20, 275--283.
\item Andrews, G.E., 1981. Ramanujan's ``lost'' notebook. I. Partial $\theta$-functions. Adv. Math. 41, 137--172.
\item Andrews, G.E., Warnaar, S.O., 2007. The product of partial theta functions. Adv. Appl. Math. 39, 116--120.
\item Chen, W.Y.C., Hou, Q.H., Mu, Y.P., 2008. Nonterminating basic hypergeometric series and the $q$-Zeilberger algorithm. Proc. Edinburgh Math. Soc. 51, 609--633.
\item Gasper, G., Rahman, M., 2004. Basic Hypergeometric Series, Encyclopedia of Mathematics and its Applications, Vol. 96, second ed. Cambridge Univ. Press, Cambridge, MA.
\item Gessel, I.M., 1995. Finding identities with the WZ method. J. Symbolic Comput. 20, 537--566.
\item Gosper, R.W., 1978. Decision procedure for indefinite hypergeometric summation. Proc. Natl. Acad. Sci. USA 75 (1), 40--42.
\item Koepf, W., 2003. Computer algebra algorithms for orthogonal polynomials and special functions. In: Koelink, E., Van Assche, W. (Eds.), Orthogonal Polynomials and Special Functions: Leuven 2002. In: Lecture Notes Math., vol. 1817. pp. 1--24.
\item Koornwinder, T.H., 1993. On Zeilberger's algorithm and its $q$-analogue. J. Comput. Appl. Math. 48, 91--111.
\item Koornwinder, T.H., 1998. Identities of nonterminating series by Zeilberger's algorithm. J. Comput. Appl. Math. 99, 449--461.
\item Koornwinder, T.H., 2005. \url{http://staff.science.uva.nl/~thk/specfun/compalg.html}.
\item Mohammed, M., Zeilberger, D., 2004. The Markov-WZ method. Electron. J. Combin. 11, R53.
\item Vidunas, R., 2002. A generalization of Kummer's identity. Rocky Mountain J. Math. 32, 919--936.
\item Warnaar, S.O., 2003. Partial theta functions. I. Beyond the lost book. Proc. London Math. Soc. 87, 363--395.
\item Wilf, H.S., Zeilberger, D., 1990a. Towards computerized proofs of identities. Bull. Amer. Math. Soc. 23 (1), 77--83.
\item Wilf, H.S., Zeilberger, D., 1990b. Rational functions certify combinatorial identities. J. Amer. Math. Soc. 3 (1), 147--158.
\item Wilf, H.S., Zeilberger, D., 1992. An algorithmic proof theory for hypergeometric (ordinary and ``$q$'') multisum/integral identities. Invent. Math. 108, 575--633.
\item Zeilberger, D., 1990. A fast algorithm for proving terminating hypergeometric identities. Discrete Math. 80, 207--211.
\item Zeilberger, D., 1994. Towards a WZ proof of Mehta's integral. SIAM J. Math. Anal. 25, 812--814.
\item Zeilberger, D., Private communication, 2008.
\end{list}\endgroup

\newpage
\section{Addendum}

 We adopt the notation of the original paper, with the assumption
$0<|q|<1$. In the examples of constructing
$q$-WZ pairs for infinite series, 
the determination of the constants requires further justification. 
The gaps can be fixed as follows.

In Example 3.1, define
\begin{align}
  F(n,k)
= 
\frac{(b,aq^n;q)_k(c/ab,cq^n;q)_\infty}
{(q,cq^n;q)_k(c/a,cq^n/b;q)_\infty}
\left(\frac{c}{ab}\right)^k,
\end{align}
where $\left|\frac{c}{ab}\right|<1$.
We need to verify that
\begin{align}
\sum_{k=-\infty}^{\infty}F(0,k)
=\sum_{k=0}^{\infty}
\frac{(a,b;q)_k(c,c/ab;q)_\infty}
{(q,c;q)_k(c/a,c/b;q)_\infty}
\left(\frac{c}{ab}\right)^k
=1.
\label{ex:3.1}
\end{align}
Replacing $a$ and $c$ with $aq^n$ and $cq^n$,
respectively,  by the $q$-WZ relation, the sum
\begin{align}
\sum_{k=-\infty}^{\infty}F(n,k)
=\sum_{k=0}^{\infty}
\frac{(b,aq^n;q)_k(c/ab,cq^n;q)_\infty}
{(q,cq^n;q)_k(c/a,cq^n/b;q)_\infty}
\left(\frac{c}{ab}\right)^k
\end{align}
 is independent of $n$,
but may depend on $a,b,c$. Hence,
\begin{align}
\sum_{k=-\infty}^{\infty}F(0,k)
=\lim_{n\to\infty}
\sum_{k=-\infty}^{\infty}F(n,k).
\end{align}
Since $F(n,k)=0$ for $k<0$, by dominated convergence
and the $q$-binomial theorem, we obtain
\begin{align}
\lim_{n\to\infty}\sum_{k=0}^{\infty}F(n,k)
&=\sum_{k=0}^{\infty}\lim_{n\to\infty}F(n,k)
\notag\\[6pt]
&=\sum_{k=0}^{\infty}
\frac{(b;q)_k(c/ab;q)_\infty}
{(q;q)_k(c/a;q)_\infty}
\left(\frac{c}{ab}\right)^k
=1.
\end{align}
This proves \eqref{ex:3.1}.

As for Example 3.2, define
\begin{align}
F(n,k)
&=
\frac{(1-aq^{n+2k})(c,d,aq^n,bq^n;q)_k}
{(1-aq^n)(q,aq/b,aq^{n+1}/c,aq^{n+1}/d;q)_k}
\notag\\[6pt]
&\qquad \times
\frac{(aq/b,aq/bcd,aq^{n+1}/c,aq^{n+1}/d;q)_\infty}
{(aq/bc,aq/bd,aq^{n+1},aq^{n+1}/cd;q)_\infty}
\left(\frac{aq}{bcd}\right)^k,
\end{align}
where $\left|\frac{aq}{bcd}\right|<1$.
We need to verify that
\begin{align}
\sum_{k=-\infty}^{\infty}F(0,k)
&=\sum_{k=0}^{\infty}
\frac{(1-aq^{2k})(a,b,c,d;q)_k}
{(1-a)(q,aq/b,aq/c,aq/d;q)_k}
\notag\\[6pt]
&\qquad \times
\frac{(aq/b,aq/c,aq/d,aq/bcd;q)_\infty}
{(aq,aq/bc,aq/bd,aq/cd;q)_\infty}
\left(\frac{aq}{bcd}\right)^k
=1.
\label{ex:3.2}
\end{align}
Replacing $a$ and $b$ with $aq^n$ and $bq^n$,
respectively, we consider
\begin{align}
\sum_{k=-\infty}^{\infty}F(n,k)
&=\sum_{k=0}^{\infty}
\frac{(1-aq^{n+2k})(c,d,aq^n,bq^n;q)_k}
{(1-aq^n)(q,aq/b,aq^{n+1}/c,aq^{n+1}/d;q)_k}
\notag\\[6pt]
&\qquad \times
\frac{(aq/b,aq/bcd,aq^{n+1}/c,aq^{n+1}/d;q)_\infty}
{(aq/bc,aq/bd,aq^{n+1},aq^{n+1}/cd;q)_\infty}
\left(\frac{aq}{bcd}\right)^k.
\end{align}
By the $q$-WZ relation, this sum is independent of $n$,
but may depend on $a,b,c,d$. Hence,
\begin{align}
\sum_{k=-\infty}^{\infty}F(0,k)
=\lim_{n\to\infty}
\sum_{k=-\infty}^{\infty}F(n,k).
\end{align}
Since $F(n,k)=0$ for $k<0$, by dominated convergence
and the $q$-Gauss summation, we obtain
\begin{align}
\lim_{n\to\infty}\sum_{k=0}^{\infty}F(n,k)
&=\sum_{k=0}^{\infty}\lim_{n\to\infty}F(n,k)
\notag\\
&=\sum_{k=0}^{\infty}
\frac{(c,d;q)_k(aq/b,aq/bcd;q)_\infty}
{(q,aq/b;q)_k(aq/bc,aq/bd;q)_\infty}
\left(\frac{aq}{bcd}\right)^k
=1.
\end{align}
This proves \eqref{ex:3.2}.

Regarding Example 3.3,  we define
\begin{align}
F(n,k)
&=
\frac{(aq^n;q)_k}{(bq^n;q)_k}
\frac{(z,b/(az),bq^n,q^{1-n}/a;q)_\infty}
{(q,b/a,azq^n,q^{1-n}/(az);q)_\infty}
z^k.
\end{align}
We need to verify that for $|z|<1$,
\begin{align}
\sum_{k=0}^{\infty}F(0,k)
&=
\sum_{k=-\infty}^{\infty}
\frac{(a;q)_k}{(b;q)_k}
\frac{(b,q/a,z,b/(az);q)_\infty}
{(q,b/a,az,q/(az);q)_\infty}
z^k
=1.
\label{ex:3.3}
\end{align}

Replacing $a$ and $b$ with $aq^n$ and $bq^n$,
respectively, we consider
\begin{align}
\sum_{k=-\infty}^{\infty}F(n,k).
\end{align}
By the $q$-WZ relation, this sum is independent
of $n$, but may depend on $a,b,z$.

Letting $b \to q$, we have for $k <0$,
\begin{align}
\lim_{b\to q}F(0,k)=0.
\end{align}
By the $q$-binomial theorem (3.9), we have
\begin{align}
\lim_{b\to q}\sum_{k=-\infty}^{\infty}F(0,k)
&=
\frac{(z;q)_\infty}{(az;q)_\infty}
\sum_{k=0}^{\infty}
\frac{(a;q)_k}{(q;q)_k}z^k=1.
\end{align}

Let
\[
H(a,b,z)=\sum_{k=-\infty}^{\infty}F(0,k).
\]
By the $q$-WZ relation, we have
\begin{align}
H(a,b,z)=H(aq^n,bq^n,z),
\qquad n\ge0.
\end{align}
Replacing $a$ with $aq^{-m}$ in the special
case $b=q$ and taking $n=m$, we obtain
\begin{align}
H(a,q^{m+1},z)
=H(aq^{-m},q,z)=1
\end{align}
for all sufficiently large integers $m$.

Since $q^{m+1}\to0$ and both sides of (3.14)
are holomorphic in $b$ near $0$ for 
$a,z$ satisfying $|q/a|<|z|<1$, the identity theorem for holomorphic functions, together with meromorphic continuation, establishes
(3.14) for general admissible parameters.
This proves \eqref{ex:3.3}.

In the case of Example 3.4, let $R(a,b,c,d,e)$ denote
the product on the right side of (3.17),
and define
\begin{align}
F(n,k)
&=
\frac{1-aq^{n+2k}}{1-aq^n}
\frac{(bq^n,cq^n,d,e;q)_k}
{(aq/b,aq/c,aq^{n+1}/d,aq^{n+1}/e;q)_k}
\notag\\
&\qquad \times
\left(\frac{a^2q}{bcde}\right)^k
\frac{1}{R(aq^n,bq^n,cq^n,d,e)}.
\end{align}
We need to verify that
\begin{align}
\sum_{k=-\infty}^{\infty}F(0,k)
&=
\frac{1}{R(a,b,c,d,e)}
\sum_{k=-\infty}^{\infty}
\frac{1-aq^{2k}}{1-a}
\notag\\[6pt]
&\qquad \times
\frac{(b,c,d,e;q)_k}
{(aq/b,aq/c,aq/d,aq/e;q)_k}
\left(\frac{a^2q}{bcde}\right)^k
=1,
\label{ex:3.4}
\end{align}
where $\left|\frac{a^2q}{bcde}\right| <1$.

Replacing $a,b,c$ with $aq^n,bq^n,cq^n$,
respectively, we consider
\begin{align}
\sum_{k=-\infty}^{\infty}F(n,k).
\end{align}
By the $q$-WZ relation, this sum is independent
of $n$, but may depend on $a,b,c,d,e$.

Letting $b\to a$, we have for $k <0$,
\[
F(0,k)=0.
\]
By the ${}_6\phi_5$ summation (3.13), we obtain
\begin{align}
\lim_{b\to a}\sum_{k=-\infty}^{\infty}F(0,k)
&=
\frac{(aq/c,aq/d,aq/e,aq/cde;q)_\infty}
{(aq,aq/cd,aq/ce,aq/de;q)_\infty}
\notag\\
&\qquad \times
\sum_{k=0}^{\infty}
\frac{(1-aq^{2k})(a,c,d,e;q)_k}
{(1-a)(q,aq/c,aq/d,aq/e;q)_k}
\left(\frac{aq}{cde}\right)^k
\notag\\[6pt]
&=1.
\end{align}

Let \[ H(a,b,c,d,e)=\sum_{k=-\infty}^{\infty}F(0,k).\]
By the same reasoning as before, we have for $n\geq 0$, 
\begin{align}
H(a,b,c,d,e)
=H(aq^n,bq^n,cq^n,d,e).
\end{align}
Since $H$ is symmetric in $b,c,d,e$,
we see that for $n\geq 0$,
\begin{align}
H(a,b,c,d,e)
=H(aq^n,b,cq^n,dq^n,e).
\end{align}

Replacing $(a,c,d)$ with
$(aq^{-m},cq^{-m},dq^{-m})$ and taking $n=m$,
we obtain from the special case that $H(a,a,c,d,e)=1$,
\begin{align}
H(a,aq^{-m},c,d,e)
&=H(aq^{-m},aq^{-m},
cq^{-m},dq^{-m},e)=1
\label{24}
\end{align}
for all sufficiently large integers $m$.

Writing $t=1/b$, then $b=aq^{-m}$ can be expressed as $t=q^m/a$.
For parameters $a,c,d,e$ satisfying
$|aq/(cde)|<1$, both sides
of (3.17) are holomorphic in $t$ near $0$.
Indeed, after the cancellation of powers of $t$,
the terms with indices $\pm k$, $k\ge1$,
are bounded on $|t|\le r$ by
\[
C_r
\left(\left|\frac{a^2q}{cde}\right|r\right)^k,
\]
where $r>0$ is sufficiently small such that
all denominators are nonzero on the disk and
$\left|a^2q/(cde)\right|r<1$.
Here $C_r>0$ is independent of $k$ and $t$.
Thus both tails converge locally uniformly.

Since $t=q^m/a\to0$, the identity theorem
for holomorphic functions and
meromorphic continuation then establishes (3.17)
for general admissible parameters.
This proves \eqref{ex:3.4}.

The last displayed equation in Example~3.2 should read
\begin{align}
&\sum_{j=0}^k
\frac{(1-aq^{2j})(a,b,c,d;q)_j}
{(1-a)(q,aq/b,aq/c,aq/d;q)_j}
\left(\frac{aq}{bcd}\right)^j
\notag\\[6pt]
&\qquad=
\frac{(aq,aq/cd;q)_\infty}
{(aq/c,aq/d;q)_\infty}
\sum_{j=0}^k
\frac{(c,d;q)_j}{(q,aq/b;q)_j}
\left(\frac{aq}{bcd}\right)^j
\notag\\[6pt]
&\qquad \qquad +
\frac{b(aq;q)_k(b,c,d;q)_{k+1}}
{(q,aq/b;q)_k(aq/c,aq/d;q)_{k+1}}
\left(\frac{aq}{bcd}\right)^{k+1}
\notag\\[6pt]
&\qquad\qquad\times
\sum_{n=0}^\infty
\frac{(aq/cd;q)_n
(aq^{k+1},bq^{k+1};q)_n}
{(b;q)_{n+1}
(aq^{k+2}/c,aq^{k+2}/d;q)_n}q^n.
\end{align}

The displayed equation immediately following (3.18) should read
\begin{align}
\lim_{b\to a}\sum_{k=-\infty}^{\infty}F(0,k)
&=
\sum_{k=0}^{\infty}
\frac{(1-aq^{2k})(a,c,d,e;q)_k}
{(1-a)(q,aq/c,aq/d,aq/e;q)_k}
\notag\\[6pt]
&\qquad\times
\frac{(aq/c,aq/d,aq/e,aq/cde;q)_\infty}
{(aq,aq/cd,aq/ce,aq/de;q)_\infty}
\left(\frac{aq}{cde}\right)^k
=1.
\end{align}

\textbf{Acknowledgments.}
We wish to thank Heng Huat Chan and Song Heng Chan for pointing out 
the incompleteness of the proofs in the examples.

\end{document}